\documentclass[runningheads]{llncs}

\usepackage{tikz}

\pgfdeclarelayer{edgelayer}
\pgfdeclarelayer{nodelayer}
\pgfsetlayers{edgelayer,nodelayer,main}

\tikzstyle{none}=[inner sep=0pt]
\definecolor{hexcolor0xf81e1c}{rgb}{0.973,0.118,0.110}
\definecolor{hexcolor0x3c00ff}{rgb}{0.235,0.000,1.000}

\tikzstyle{whitevertex}=[circle,fill=white,draw=black, scale = 0.5]
\tikzstyle{redvertex}=[circle,fill=hexcolor0xf81e1c,draw=black, scale = 0.5]
\tikzstyle{bluevertex}=[circle,fill=hexcolor0x3c00ff,draw=black, scale = 0.5]
\tikzstyle{greenvertex}=[circle,fill=green,draw=black, scale=0.5]
\tikzstyle{purplevertex}=[circle,fill=magenta,draw=black, scale=0.5]
\tikzstyle{grayvertex}=[circle,fill=white,draw=gray, scale=0.5]
\tikzstyle{blackvertex}=[circle,fill=black,draw=black, scale=0.5]
\tikzstyle{doublearc}=[black, <->]

\tikzstyle{textbox}=[rectangle,fill=none,draw=none]
\tikzstyle{box}=[rectangle,fill=none,draw=black]

\tikzstyle{arc}=[black, ->]
\tikzstyle{grayarc}=[gray, ->]
\tikzstyle{bluearc}=[blue, ->]
\tikzstyle{grayedge}=[draw=gray]
\tikzstyle{blueedge}=[draw=blue]
\tikzstyle{rededge}=[draw=red]
\tikzstyle{edge}=[draw=black]

\tikzstyle{vertex}=[circle, ,fill=white,draw=black, scale=0.5]

\tikzstyle{10circle}=[circle, scale=10.0,draw=black]
\tikzstyle{10oval}=[ellipse, scale=10.0,draw=black]

\usepackage{amssymb}
\usepackage{graphicx}
\usepackage{url}
\urldef{\mailsa}\path|huangj@uvic.ca|
\urldef{\mailsb}\path|fayye@uvic.ca|

\newcommand{\keywords}[1]{\par\addvspace\baselineskip
\noindent\keywordname\enspace\ignorespaces#1}

\usepackage{amsmath, amssymb}

\begin{document}

\mainmatter  


\title{Semi-strict chordality of digraphs}


\titlerunning{Semi-strict chordality}

%
%
\author{Jing Huang\inst{1} \and Ying Ying Ye\inst{2}}

\authorrunning{J. Huang and Y.Y. Ye }


\institute{Department of Mathematics and Statistics, University of Victoria,\\
                  Victoria, B.C., Canada,\ V8W\ 2Y2\\
\mailsa
\and Department of Mathematics and Statistics, University of Victoria,\\
                  Victoria, B.C., Canada,\ V8W\ 2Y2\\
\mailsb}


%
%

\toctitle{Lecture Notes in Computer Science} 
\tocauthor{Jing Huang and Ying Ying Ye} 

\maketitle

\begin{abstract}
Chordal graphs are important in algorithmic graph theory. Chordal digraphs are 
a digraph analogue of chordal graphs and have been a subject of active studies 
recently. Unlike chordal graphs, chordal digraphs lack many structural properties
such as forbidden subdigraph or representation characterizations. In this paper we 
introduce the notion of semi-strict chordal digraphs which form a class strictly 
between chordal digraphs and chordal graphs. Semi-strict chordal digraphs have 
rich structural properties. We characterize semi-strict chordal digraphs in terms
of knotting graphs, a notion alalogous to the one introduced by Gallai for the 
study of comparability graphs. We also give forbidden subdigraph characterizations
of semi-strict chordal digraphs within the cases of locally semicomplete digraphs 
and weakly quasi-transitive digraphs.
\end{abstract}

\section{Introduction}

Amongst special classes of graphs, chordal graphs are perhaps the most important in
the structural graph theory, cf. \cite{golumbic,rs}. They are the graphs which do 
not contain an induced cycle of length four or more. They can be equivalently 
defined in terms of existence of simplicial vertices. A {\em simplicial} vertex in 
a graph is a vertex whose neighbourhood induces a clique. A graph is {\em chordal} 
if and only if its every induced subgraph contains a simplicial vertex. 

The digraph analogue of chordal graphs has been introduced in \cite{HR} and studied
in \cite{hy,kleitman,meister,ye}. 
A vertex $v$ in a digraph $D$ is called {\em di-simplicial} if for any 
$u \in N^-(v)$ and $w \in N^+(v)$ with $u \neq w$, $uw$ is an arc in $D$. 
A digraph is called {\em chordal} if its every induced subdigraph contains 
a di-simplicial vertex. Chordal graphs are essentially the symmetric chordal 
digraphs and form a proper subclass of chordal digraphs. 

Unlike chordal graphs, which admit nice characterizations (cf. \cite{golumbic}), 
the structure of chordal 
digraphs is complex. In particular, there is no known forbidden subdigraph 
characterization for chordal digraphs in general. Meister and Telle \cite{meister}
characterized semicomplete chordal digraphs by forbidden subdigraphs. 
Locally semicomplete digraphs and weakly quasi-transitive digraphs are two 
classes of digraphs, both larger than the class of semicomplete digraphs 
\cite{bang,hy}. Generalizing the result of Meister and Telle in \cite{meister},
the authors gave a forbidden subdigraph characterizations of chordal digraphs 
within these two classes of digraphs \cite{hy}. However, it remains an open problem
to find a forbidden subdigraph characterization for chordal digraphs in general

Several variants of chordal digraphs have been introduced in \cite{2h,ye}.
One success along this line of study is the introduction of the notion of strict
chordal digraphs \cite{2h}. A vertex $v$ in a digraph $D$ is called 
{\em strict di-simplicial} if for any $u, w \in N^-(v) \cup N^+(v)$ with $u \neq w$,
both $uw, wu$ are arcs in $D$. A digraph is called {\em strict chordal} if its
every induced subdigraph contains a strict di-simplicial vertex. A forbidden
subdigraph characterization of strict chordal digraphs has been obtained in 
\cite{2h}. The paper \cite{mckee} provides an interesting point of view of 
the structure of strict strict chordal digraphs. 

Although strict chordal digraphs contain all chordal graphs, they are a much 
restricted subclass of chordal digraphs. In this paper, we propose a notion of 
digraphs that is between chordal digraphs and strict chordal digraphs. 
A vertex $v$ in a digraph $D$ is called {\em semi-strict di-simplicial} if 
for any $u \in N^-(v)$ and $w \in N^+(v)$ with $u \neq w$, both $uw, wu$ are arcs 
in $D$. A digraph is called {\em semi-strict chordal} if its every induced 
subdigraph contains a semi-strict di-simplicial vertex. It follows from definition 
that every semi-strict chordal digraph $D$ has a vertex ordering 
$v_1, v_2, \dots, v_n$ such that $v_i$ is a semi-strict di-simplicial vertex
in the subdigraph of $D$ induced by $v_i, v_{i+1}, \dots, v_n$ for each $i \geq 1$.
Such an ordering is called a {\em semi-strict perfect elimination ordering} of $D$.
It is clear that if a digraph has a semi-strict perfect elimination ordering then
it is semi-strict chordal. Hence a digraph is semi-strict chordal if and only if
it has a semi-strict perfect elimination ordering. A semi-strict perfect 
elimination ordering of a digraph can be obtained (if one exists) by successively
finding a semi-strict di-simplicial vertex, which can be done in polynomial time.
Therefore semi-strict chordal digraphs can be recognized in polynomial time.   

We characterize semi-strict chordal digraphs in terms of their knotting graphs, 
which are auxiliary graphs defined for digraphs (see Section \ref{two}). 
The knotting graphs for digraphs are analogous to the knotting graphs
for graphs introduced by Gallai \cite{gallai} of the study of comparability graphs
(i.e., the graphs which have transitive orientations). Knotting graphs for graphs
are a useful concept as Gallai proved that a graph is a comparability graph if and
only if its knotting graph is bipartite. We demonstrate the usefulness of the 
knotting graphs for digraphs by showing that a digraph is semi-strict chordal if
and only if every induced subgraph of its knotting graph has a subset of vertices
whose degrees are either zero or one.

Semi-strict chordal digraphs are still general enough. It appears challenging to 
find a forbidden subdigraph characterization for the class. The characterization of
semi-strict chordal digraphs in terms of knotting graphs does not seem to yield
a forbidden subdigraph characterization of the class. So we turn our attention
to locally semicomplete digraphs and weakly quasi-transitive digraphs, for which
forbidden subdigraph characterizations of chordal digraphs are known \cite{hy}.
We will show that semi-strict chordal digraphs can be characterized by forbidden
subdigraphs within these two classes of digraphs.

All digraphs considered in this paper do not contain loops or multiple arcs but may 
contain {\em digons} (i.e., pairs of arcs joining vertices in opposite directions).
If an arc is contained in a digon then it is called a {\em symmetric arc}. A digraph
which does not contain any symmetric arc is called an {\em oriented graph}. 
A digraph which contains only symmetric arcs is called a {\em symmetric digraph}. 
Graphs may be viewed as symmetric digraphs.

Two vertices in a digraph $D$ are {\em adjacent} and referred to as {\em neighbours}
of each other if there is at least one arc between them. We say that $u$ is 
an {\em in-neighbour} of $v$ or $v$ an {\em out-neighbour} of $u$ if $uv$ is an arc 
in $D$ (symmetric or not). The set of all in-neighbours of a vertex $v$ is denoted 
by $N^-(v)$ and the set of all out-neighbours of $v$ is denoted by $N^+(v)$. 

Let $v$ be a vertex and $u, w$ be neighbours of $v$ in a digraph $D$. 
Then $u, w$ are {\em synchronous neighbours} of $v$ if $u, w$ are both in $N^-(v) \setminus N^+(v)$, or in $N^+(v) \setminus N^-(v)$, or in $N^-(v) \cap N^+(v)$;
otherwise they are called {\em asynchronous} neighbours of $v$.
We use $S(D)$ to denote the spanning subdigraph of $D$ whose arc set consists of all
symmetric arcs in $D$.

A digraph $D$ is {\em semicomplete} if any two vertices are adjacent.
A digraph $D$ is called {\em locally semicomplete} if for every vertex $v$, $N^-(v)$
and $N^+(v)$ each induces a semicomplete subdigraph in $D$. 
A digraph $D$ is called {\em weakly quasi-transitive} if for each vertex $v$ of $D$, any two asynchronous neighbours of $v$ are adjacent.

Weakly quasi-transitive digraphs generalize extended semicomplete digraphs as well
as quasi-transitive digraphs which are studied in \cite{BH}. 
A digraph is called an {\em extended semicomplete} digraph if it is obtained from
a semicomplete digraph by substituting an independent set for each vertex of 
the semicomplete digraph. A digraph is {\em qausi-transitive} if for any 
three vertices $u, v, w$ if $u \in N^-(v)$ and $w \in N^+(v)$ then $u$ and $w$ are
adjacent. As shown in \cite{BH}, quasi-transitive digraphs 
can be constructed recursively from transitive oriented graphs and semicomplete 
digraphs by substitutions. It turns out weakly quasi-transitive digraphs can also 
be constructed in a similar way as shown in \cite{hy}. 

\begin{theorem}[\cite{hy}]\label{WQTDStructure}
Let $\cal{W}$ be the class of digraphs defined recursively as follows:
\begin{tabbing}
\= $\bullet$ each transitive oriented graph is in $\cal{W}$;\\
\> $\bullet$ each semicomplete digraph is in $\cal{W}$;\\
\> $\bullet$ each symmetric digraph is in $\cal{W}$;\\
\> $\bullet$ if $D$ is in $\cal{W}$, then any digraph obtained from $D$ by
          substituting digraphs in $\cal{W}$\\
\> \hspace{2mm} for the vertices of $D$ is in $\cal{W}$.
\end{tabbing}
Then $\cal{W}$ is precisely the class of weakly quasi-transitive digraphs.\qed
\end{theorem}

The paper is organized as follows. In Section \ref{two}, we introduce the concept
of knotting graphs for digraphs and use it to characterize semi-strict chordal
digraphs. We will give forbidden subdigraph characterizations of weakly 
quasi-transitive semi-strict chordal digraphs and locally semicomplete 
semi-strict chordal digraphs in Sections \ref{three} and \ref{four} respectively.
Finally, in Section \ref{five}, we finish with some concluding remarks and open 
problems.

\section{Knotting graphs}
\label{two}

In his study of comparability graphs (i.e., transitively orientable graphs), 
Gallai \cite{gallai} introduced the concept of knotting graphs. It turns out this 
is a very useful concept. Each graph gives rise to a unique knotting graph and 
a graph is a comparability graph if and only if its corresponding knotting graph is
bipartite, cf. \cite{gallai}. As knotting graphs can be constructed in polynomial 
time, this gives a way to recognize comparability graphs in polynomial time.

In this section, we introduce the knotting graphs for digraphs analogous to 
Gallai's knotting graphs for graphs. As a result, semi-strict chordal digraphs
can be also be recognized via their knotting graphs. 

Let $D$ be a digraph and $v$ be a vertex of $D$. Denote by $E_v$ the set of arcs
incident with $v$ (i.e., have $v$ as an endvertex). An arc $e$ in $E_v$ is called
{\em out-going} if $v$ is the tail of $e$; otherwise it is called {\em in-coming}. 
Two arcs $e, f \in E_v$ are called {\em knotted}, denoted by $e \sim f$, if there 
exists a sequence of arcs
\[e_1, e_2, \dots, e_k \]
where $e_1 = e$ and $e_k = f$ such that, for each $i$,
\begin{tabbing}
\= $\bullet$ $e_i \in E_v$, \\
\> $\bullet$ $e_i$ and $e_{i+1}$ have different endvertices and are not both 
    in-coming or both \\
\> \hspace{2mm} out-going, and\\
\> $\bullet$ the two endvertices of $e_i$ and $e_{i+1}$ distinct from $v$ are not 
   connected by a digon. 
\end{tabbing}
It is easy to verify that $\sim$ is an equivalence relation on $E_v$ and hence 
partitions $E_v$ into equivalence classes. 

Denote by $v^1, v^2, \dots, v^{\ell_v}$ the equivalence classes of $E_v$ for each
vertex $v$ of $D$. In the case when $v$ is an isolated vertex (i.e., 
$E_v = \emptyset$, we let $\ell_v = 1$ and $v^1 = \emptyset$. 

Define the {\em knotting graph} $K_D$ of $D$ as follows:
The vertex set of $K_D$ consists of $v^1, v^2, \dots, v^{\ell_v}$ for all vertices
$v$ of $D$, and the edge set of $K_D$ consists of $u^iv^j$ such that $u \neq v$
and $u^i \cap v^j \neq \emptyset$ for all $1 \leq i \leq \ell_u$ and 
$1 \leq j \leq \ell_v$. See Fig. \ref{example1} for an example of the knotting
graph of a digraph.

\begin{center} 
	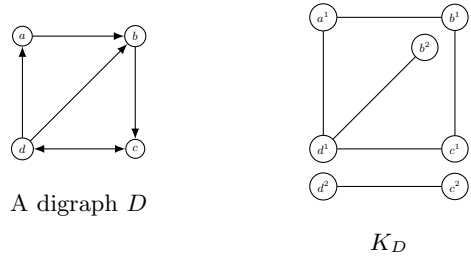
\begin{figure}[htb]
		\center
		\begin{tikzpicture}[>=latex]
			\begin{pgfonlayer}{nodelayer}
				
				\node [style=vertex] (a) at (0,0) {$a$};
				\node [style=vertex] (b) at (1.5,0) {$b$};
				\node [style=vertex] (c) at (1.5,-1.5) {$c$};
				\node [style=vertex] (d) at (0,-1.5) {$d$};
				
				\node [style=textbox] () at (0.75,-2.25) {A digraph $D$};

				\draw[style=arc] (a) to (b);
				\draw[style=arc] (b) to (c);
				\draw[style=arc] (d) to (a);
				\draw[style=arc] (d) to (b);
				\draw[style=doublearc] (d) to (c);

				\node [style=vertex] (a1) at (4, 0.25) {$a^1$};
				\node [style=vertex] (b1) at (5.75, 0.25) {$b^1$};
				\node [style=vertex] (b2) at (5.35, -0.15) {$b^2$};
				\node [style=vertex] (c1) at (5.75, -1.5) {$c^1$};
				\node [style=vertex] (c2) at (5.75, -2) {$c^2$};
				\node [style=vertex] (d1) at (4, -1.5) {$d^1$};
				\node [style=vertex] (d2) at (4, -2) {$d^2$};
				
				\node [style=textbox] () at (4.875,-2.75) {$K_D$};
				
				\draw[style=edge] (a1) to (b1);
				\draw[style=edge] (c1) to (b1);
				\draw[style=edge] (d1) to (c1);
				\draw[style=edge] (a1) to (d1);
				\draw[style=edge] (b2) to (d1);
				\draw[style=edge] (c2) to (d2);
				
			\end{pgfonlayer}
		\end{tikzpicture}
\caption{\label{example1}A digraph $D$ and its knotting graph $K_D$ where
$a^1 = \{ab,da\}, b^1 = \{ab,bc\}, b^2 = \{db\}, c^1 = \{bc,cd\}, c^2 = \{dc\},
d^1 = \{da,cd,db\}, \mbox{and} d^2 = \{dc\}$.}
	\end{figure}
\end{center}

Note that each vertex $v$ of $D$ splits into $\ell_v$ vertices 
$v^1, v^2, \dots, v^{\ell_v}$ in $K_D$.  Each arc in $E_v$ uniquely corresponds to
an edge incident with $v^j$ for some $1 \leq j \leq \ell_v$. In general, each arc 
of $D$ with endvertices $u, v$ uniquely corresponds to an edge $u^iv^j$ of $K_D$ 
for some $1 \leq i \leq \ell_u$ and $1 \leq j \leq \ell_v$. This correspondence is 
one-to-one between the arc set of $D$ and the edge set of $K_D$. Moreover, for each
vertex $v$ of $D$, the knotting graph of $D-v$ is the one obtained from $K_D$ by 
deleting the vertices $v^1, v^2, \dots, v^{\ell_v}$. Thus, for each induced 
subdigraph $H$ of $D$, the knotting graph $K_H$ of $H$ can be computed efficiently 
from $K_D$.     

\begin{lemma} \label{one}
Let $D$ be a digraph and $K_D$ be the knotting graph of $D$. Then a vertex $v$ of 
$D$ is semi-strict di-simplicial if and only if its splitting vertices 
$v^1, v^2, \dots, v^{\ell_v}$ all have degree zero or one in $K_D$.
\end{lemma}

\begin{proof}
Suppose that $v$ is a semi-strict di-simplicial vertex. If $v$ is isolated, then
its unique splitting vertex $v^1$ has degree zero. So assume that $v$ is not 
isolated. By definition, for any $u \in N^-(v)$ and $w \in N^-(v)$ with $u \neq w$,
$u$ and $w$ are connected by a digon. It follows that each arc in $E_v$ is knotted 
only with itself. So $v$ splits into $|N^-(v)|+|N^+(v)|$ vertices in $K_D$, each 
has degree one.

On the other hand, suppose that $v$ is not a semi-strict di-simplicial vertex. 
Then for some $u \in N^-(v)$ and $w \in N^-(v)$, $u$ and $w$ are not connected by
a digon. It follows that the two arcs $uv$ and $vw$ are knotted, that is, they are 
in the same equivalence class $v^j$ for some $1 \leq j \leq \ell_v$. Hence $uv, vw$
correspond to two edges incident with $v^j$ in $K_D$, showing that the degree of 
$v^j$ is at least two.  
\qed
\end{proof}

A {\em splitting group} in $K_D$ consists of all splitting vertices 
$v^1, v^2, \dots, v^{\ell_v}$ for some vertex $v$ of $D$. For an induced subgraph 
$S$ of $K_D$, a {\em splitting group} in $S$ consists of all vertices of $S$ which 
are in some splitting group of $K_D$.

\begin{theorem} 
\label{degree-condition}
A digraph $D$ is semi-strict chordal if and only if every induced subgraph of $K_D$
has a splitting group whose vertices each has degree zero or one.
\end{theorem}

\begin{proof}
Suppose that $D$ is semi-strict chordal. Let $S$ be a subgraph of $K_D$ and let
$H$ be the subdigraph of $D$ induced by the vertices which have splitting vertices
in $S$. Since $D$ is semi-strict chordal, $H$ has a semi-strict di-simplicial 
vertex $v$. It follows from Lemma \ref{one} that the splitting vertices of $v$
in $S$ each has degree zero or one. 

Conversely, suppose that every induced subgraph of $K_D$ has a splitting group 
whose vertices each has degree zero or one. Let $H$ be an induced subdigraph of 
$D$. Then $K_H$ is an induced subgraph of $K_D$ and has a splitting group whose
vertices have degree zero or one. By Lemma \ref{one} the vertex corresponding to 
the splitting group is a semi-strict di-simplicial vertex in $H$. Therefore
$D$ is a semi-strict chordal digraph.   
\qed
\end{proof}

For the knotting graph $K_D$ in Fig. \ref{example1}, each splitting group has a 
vertex of degree 2 or more, so by Theorem \ref{degree-condition} the corresponding
digraph $D$ is not semi-strict chordal. In Fig. \ref{example2}, the knotting graph 
$K_D$ satisfies the property that every induced subgraph has a splitting group 
whose vertices each has degree zero or one, hence according to 
Theorem \ref{degree-condition} the digraph $D$ is semi-strict chordal. 
    
\begin{center} 
	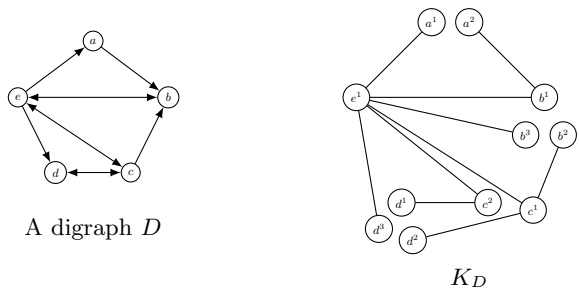
\begin{figure}[htb]
		\center
		\begin{tikzpicture}[>=latex]
			\begin{pgfonlayer}{nodelayer}
				
				\node [style=vertex] (a) at (0,0.75) {$a$};
				\node [style=vertex] (b) at (1,0) {$b$};
				\node [style=vertex] (c) at (0.5,-1) {$c$};
				\node [style=vertex] (d) at (-0.5,-1) {$d$};
				\node [style=vertex] (e) at (-1,0) {$e$};
				
				\node [style=textbox] () at (0,-1.75) {A digraph $D$};

				\draw[style=arc] (a) to (b);
				\draw[style=arc] (e) to (a);
				\draw[style=arc] (c) to (b);
				\draw[style=doublearc] (b) to (e);
				\draw[style=doublearc] (d) to (c);
				\draw[style=doublearc] (e) to (c);
				\draw[style=arc] (e) to (d);
				\node [style=vertex] (a1) at (4.5, 1) {$a^1$};	
				\node [style=vertex] (a2) at (5, 1) {$a^2$};
				\node [style=vertex] (b1) at (6, 0) {$b^1$};	
				\node [style=vertex] (b2) at (6.25, -0.5) {$b^2$};		
				\node [style=vertex] (b3) at (5.75, -0.5) {$b^3$};	
				\node [style=vertex] (c2) at (5.25, -1.4) {$c^2$};	
				\node [style=vertex] (c1) at (5.85, -1.5) {$c^1$};		
				\node [style=vertex] (d1) at (4.1, -1.4) {$d^1$};
				\node [style=vertex] (d2) at (4.25, -1.9) {$d^2$};
				\node [style=vertex] (d3) at (3.8, -1.75) {$d^3$};
				\node [style=vertex] (e1) at (3.5, 0) {$e^1$};
				
				\node [style=textbox] () at (5,-2.4) {$K_D$};
				
				\draw[style=edge] (a1) to (e1);
				\draw[style=edge] (a2) to (b1);
				\draw[style=edge] (b1) to (e1);
				\draw[style=edge] (b2) to (c1);
				\draw[style=edge] (b3) to (e1);
				\draw[style=edge] (c1) to (e1);
				\draw[style=edge] (c2) to (e1);
				\draw[style=edge] (c2) to (d1);
				\draw[style=edge] (c1) to (d2);
				\draw[style=edge] (d3) to (e1);
				
			\end{pgfonlayer}
		\end{tikzpicture}
\caption{\label{example2}A semi-strict chordal digraph $D$ and its knotting graph 
$K_D$ where $a^1 = \{ea\}, a^2 = \{ab\}, b^1 = \{ab,bc\}, b^2 = \{cb\},
b^3 = \{eb\}, c^1 = \{cb,dc,ce\}, c^2 = \{ec,cd\}, d^1 = \{ec,cd\}, 
d^2 = \{dc,cd,db\}, d^3 = \{ed\}, \mbox{and} 
e^1 = \{ea,eb,be,ec,ce,ed\}$.}
	\end{figure}
\end{center}

\section{Weakly quasi-transitive semi-strict chordal digraphs}
\label{three}

The purpose of this section is to give a forbidden subdigraph characterization
of weakly quasi-transitive semi-strict chordal digraphs.

Let $D$ be a digraph and $C: v_1v_2 \dots v_kv_1$ be a directed cycle in $D$, We 
call $C$ {\em induced} if there is no arc between $v_i$ and $v_j$ for all $i,j$ 
with $|i-j| \notin \{1,k-1\}$.

\begin{lemma}
	If $D$ is a semi-strict chordal digraph, then the underlying graph of $S(D)$ is chordal, and every semi-strict di-simplicial vertex of $D$ is a semi-strict di-simplicial vertex of $S(D)$.
\end{lemma}
\begin{proof}
	Suppose that $D$ is semi-strict chordal with a perfect elimination ordering $\prec$.
	If the underlying graph of $S(D)$ is not a chordal graph, then $S(D)$ contains a chordless cycle $C = u_1u_2\dots u_k$ of length at least $4$.
	Without loss of generality, assume that $u_1 \prec u_2$ and $u_1 \prec u_k$. 
	Since $u_1u_2, u_1u_k$ both are symmetric arcs and $D$ is semi-strict chordal, $u_2u_k$ is a symmetric arc, which contradicts with the assumption that $C$ is a chordless cycle.
	Hence, the underlying graph of $S(D)$ is a chordal graph.

	Suppose that $v$ is a semi-strict di-simplicial vertex of $D$.
	If $v$ is not a semi-strict di-simplicial vertex of $S(D)$, then there exist vertices $u$ and $w$ such that both $uv$ and $vw$ are symmetric arcs but $uw$ is not a symmetric arc in $D$, which contradicts the fact that $v$ is a semi-strict di-simplicial vertex of $D$.
\qed
\end{proof}

\begin{lemma}\label{dicycle}
	If $D$ is a semi-strict chordal digraph, then $D$ does not contain an induced directed cycle consisting of non-symmetric arcs and $S(D)$ does not contain an induced directed cycle of length $\geq 4$.
\end{lemma}

\begin{proof}
Suppose that $C$ is either an induced directed cycle in $D$ consisting of non-symmetric arcs or an induced directed cycle of length $\geq 4$ in $S(D)$.
Then the subdigraph of $D$ that induced by the vertices of $C$ has no semi-strict di-simplicial vertex and hence is not a semi-strict chordal digraph.
Therefore, $D$ is not semi-strict chordal. 
\qed 
\end{proof}

In Fig.~\ref{Forbidden1}, the adjacency between two vertices $u$ and $v$ is depicted
by a solid line that joining $u$ and $v$. Symmetric arcs between $u$ and $v$ is 
depicted by a solid line with arrows at both ends. Any non-symmetric arc from $u$
to $v$ is depicted by a solid arrowed line from $u$ to $v$.
Any non-symmetric arc between $u$ and $v$ with its direction unspecified is depicted
by a solid line with short bar in the middle (the arc between $u$ and $v$ may go 
from $u$ to $v$ or go from $v$ to $u$).

\begin{lemma}\label{noOrientedCycle}
	Let $D$ be a weakly quasi-transitive digraph that does not contain any of the digraphs in Fig.~\ref{Forbidden1} as an induced subdigraphs and $S(D)$ does not contain an induced directed cycle of length $\geq 4$. 
	Then $D$ does not contain a directed cycle consisting of non-symmetric arcs. 
\end{lemma}
\begin{proof}
Suppose that $D$ contains such cycle and let $C = v_1v_2\dots v_kv_1$ be the one with minimum length $k$.
Clearly, $k \geq 4$ as otherwise $C$ is Fig.~\ref{Forbidden1}(d).
Since $v_1$ and $v_3$ are asynchronous neighbours of $v_2$ and $D$ is a weakly quasi-transitive digraph, there is an arc between $v_1$ and $v_3$.
However, $v_1$ and $v_3$ are not joining by a non-symmetric arc, as otherwise there exists a directed cycle consisting of non-symmetric arcs with length smaller than $C$, which contradicts with the assumption.
Hence, $v_1v_3$ is a symmetric arc.
Similarly, $v_2$ and $v_4$ are joining by a symmetric arc.
Since $v_1$ and $v_4$ are asynchronous neighbours of $v_2$, they are adjacent. 
Then no matter what kind of arc is between $v_1$ and $v_4$, the subdigraph that induced by $v_1,v_2,v_3,v_4$ is Fig.~\ref{Forbidden1}(a), a contradiction.
Therefore, no such induced directed cycle with non-symmetric arcs exists.
\qed
\end{proof}

\begin{center} 
	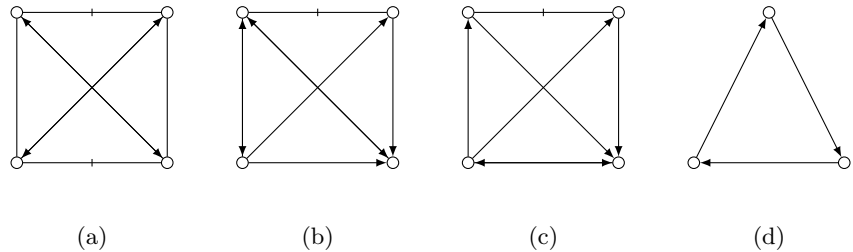
\begin{figure}[htb]
		\center
		\begin{tikzpicture}[>=latex]
		\begin{pgfonlayer}{nodelayer}
		\node [style=vertex] (1) at (0,0) {};
		\node [style=vertex] (2) at (2,0) {};
		\node [style=vertex] (3) at (2,-2) {};
		\node [style=vertex] (4) at (0,-2) {};
		
		\draw (1) to (2);
		\draw (2) to (3);
		\draw (3) to (4);
		\draw (4) to (1);
		\draw (1,-0.05) -- (1,0.05);
		\draw (1,-1.95) -- (1,-2.05);
		\draw[style=arc] (1) to (3);
		\draw[style=arc] (3) to (1);
		\draw[style=arc] (2) to (4);
		\draw[style=arc] (4) to (2);
		\node [style=textbox] at (1, -3) {(a)};

		\node [style=vertex] (1) at (3,0) {};
		\node [style=vertex] (2) at (5,0) {};
		\node [style=vertex] (3) at (5,-2) {};
		\node [style=vertex] (4) at (3,-2) {};
		
		\draw (1) to (2);
		\draw (4,-0.05) -- (4,0.05);
		\draw[style=arc] (2) to (3);
		\draw[style=arc] (4) to (3);
		\draw[style=arc] (1) to (4);
		\draw[style=arc] (4) to (1);		
		\draw[style=arc] (1) to (3);
		\draw[style=arc] (3) to (1);
		\draw[style=arc] (4) to (2);
		\node [style=textbox] at (4, -3) {(b)};
		
		\node [style=vertex] (1) at (6,0) {};
		\node [style=vertex] (2) at (8,0) {};
		\node [style=vertex] (3) at (8,-2) {};
		\node [style=vertex] (4) at (6,-2) {};
		
		\draw (1) to (2);
		\draw (7,-0.05) -- (7,0.05);
		\draw[style=arc] (2) to (3);
		\draw[style=arc] (4) to (3);
		\draw[style=arc] (3) to (4);
		\draw[style=arc] (4) to (1);
		\draw[style=arc] (1) to (3);
		\draw[style=arc] (4) to (2);
		\node [style=textbox] at (7, -3) {(c)};
		
		\node [style=vertex] (1) at (10,0) {};
		\node [style=vertex] (2) at (11,-2) {};
		\node [style=vertex] (3) at (9,-2) {};
		
		\draw[style=arc] (1) to (2);
		\draw[style=arc] (2) to (3);
		\draw[style=arc] (3) to (1);
		
		\node [style=textbox] at (10, -3) {(d)};
		
		\end{pgfonlayer}
		\end{tikzpicture}
		\caption{\label{Forbidden1}Digraphs which are not semi-strict chordal.}
	\end{figure}
\end{center}

\begin{lemma} \label{necessity} 
	If $D$ is a semi-strict chordal digraph, then $D$ does not contain any of the digraphs in Fig.~\ref{Forbidden1} as an induced subdigraph. 
\end{lemma}
\begin{proof}
Suppose that $D'$ is any of the digraphs in Fig.~\ref{Forbidden1}.
Then 
the subdigraph of $D$ induced by the vertices of $D'$ has no semi-strict di-simplicial vertex 
and hence is not a semi-strict chordal digraph. Therefore $D$ is not a semi-strict chordal digraph.  
\qed
\end{proof}

Suppose that $v$ is not a semi-strict di-simplicial vertex of $D$ but is a semi-strict di-simplicial vertex of $S(D)$.
Then there exist $u \in N^-(v)$ and $w \in N^+(v)$ such that $uw$ is not a symmetric arc of $D$ ($u$ and $w$ might not adjacent, or adjacent by a non-symmetric arc).
Such ordered triple $(u,v,w)$ of vertices is called \textit{violating triple for $v$}.
We remark that a violating triple $(u,v,w)$ exists only if $v$ is a semi-strict di-simplicial vertex of $S(D)$ and it certifies that $v$ is not a semi-strict di-simplicial vertex of $D$.
We call $v$ \textit{type 1} if for every violating triple $(u,v,w)$ both $uv$ and $vw$ are non-symmetric and \textit{type 2} otherwise.

\begin{lemma}\label{canonical_WQTD}
	Let $D$ be a weakly quasi-transitive digraph such that $S(D)$ is semi-strict chordal and $D$ does not contain any digraph in Fig.~\ref{Forbidden1} as an induced subdigraph.
	Suppose that $(u,v,w)$ is a violating triple.
	Then either $u$ or $w$ is a semi-strict di-simplicial vertex of $S(D)$, and the arc between that semi-strict di-simplicial vertex and $v$ is non-symmetric.  
\end{lemma} 
\begin{proof}
Assume that $uv$ is a non-symmetric arc (the case for $vw$ is a non-symmetric arc can be proved similarly).
Clearly, the arc between $u$ and $w$ is non-symmetric.
Consider $S(D)-(N^-[w] \cap N^+[w])$ where $N^-[w] = N^-(w) \cup \{w\}$ and  $N^+[w] = N^+(w) \cup \{w\}$.
Since $u$ and $w$ are joined by a non-symmetric arc, $u$ is a vertex of $S(D) - (N^-[w]\cap N^+[w])$.
Since $S(D)$ is semi-strict chordal, $S(D) - (N^-[w]\cap N^+[w])$ is also semi-strict chordal.
Thus each component of $S(D) - (N^-[w]\cap N^+[w])$ contains a semi-strict di-simplicial vertex. 
Let $u'$ be a semi-strict di-simplicial vertex of the component of $S(D)- (N^-[w]\cap N^+[w])$ that contains $u$. 
If $u = u'$, then we are done.
Otherwise, let $u_1u_2\dots u_k$ where $u_1=u$ and $u_k = u'$ be a directed path in $S(D) - (N^-[w]\cap N^+[w])$.
Thus, there exists vertex $u_2$ such that $uu_2$ is a symmetric arc.
Moreover, since $u_2 \in S(D) - (N^-[w]\cap N^+[w])$, and since $u_2$ and $w$ are asynchronous neighbours of $u$, they are joining by a non-symmetric arc.
Suppose that $vw$ is a symmetric arc, then since $v$ and $u_2$ are adjacent, $u,v,w,u_2$ induce Fig.~\ref{Forbidden1}(a), a contradiction.
Therefore, $vw$ is a non-symmetric arc.
Now consider $S(D)-(N^-[u] \cap N^+[u])$.
Since the arc between $u$ and $w$ is not symmetric, we know that $w \in S(D)-(N^-[u] \cap N^+[u])$, and so the component of $S(D) - (N^-[u] \cap N^+[u])$ that contains $w$ has a semi-strict di-simplicial vertex. 
Let $w'$ be such semi-strict di-simplicial vertex. 
If $w=w'$, then our statement is proved.
Hence, assume that $w_1w_2\dots w_l$ be a directed path in $S(D)-(N^-[u] \cap N^-[u])$ where $w=w_1$ and $w'=w_l$.
Then $ww_2$ is a symmetric arc.
Since $u$ and $w_2$ are asynchronous neighbours of $w$, they are adjacent by a non-symmetric arc.
Moreover, $u_2$ and $w_2$ are adjacent because they are asynchronous neighbours of $w$.
But then the subdigraph that induced by $u,u_2,w,w_2$ is Fig.~\ref{Forbidden1}(a), which is a contradiction.
Therefore, either $u$ or $w$ is a semi-strict di-simplicial vertex of $S(D)$, and it is joining with $v$ by a non-symmetric arc.
\qed
\end{proof}

We note that the above lemma works for any violating triple $(u,v,w)$.

\begin{theorem}\label{WQTDMain}
	A weakly quasi-transitive digraph $D$ is semi-strict chordal if and only if $S(D)$ is semi-strict chordal and $D$ does not contain any digraph in Fig.~\ref{Forbidden1} as an induced subdigraph. 
\end{theorem}
\begin{proof}
The necessity follows from Lemma \ref{dicycle} and Lemma \ref{necessity}.
For the other direction assume that $S(D)$ is semi-strict chordal and $D$ does not contain any digraph in Fig.~\ref{Forbidden1} as an induced subdigraph.
To prove $D$ is semi-strict chordal it suffices to show that $D$ has a semi-strict di-simplicial vertex.
Since $S(D)$ is semi-strict chordal, it has semi-strict di-simplicial vertices. 
If any semi-strict di-simplicial vertex of $S(D)$ is a semi-strict di-simplicial vertex of $D$, then we are done. 
Hence, we assume that none of the semi-strict di-simplicial vertices of $S(D)$ is a semi-strict di-simplicial vertex of $D$.

First suppose that $S(D)$ has semi-strict di-simplicial vertices of type 1. 
Let $v$ be such a vertex.
Then there exists a violating triple $(u,v,w)$ for $v$, and by Lemma \ref{canonical_WQTD} either $u$ or $w$ is a semi-strict di-simplicial vertex of $S(D)$. 
Know that both $uv$ and $vw$ are non-symmetric arcs, and so $u$ and $w$ are adjacent because they are asynchronous neighbours of $v$.
But $wu$ is not a non-symmetric arc because otherwise $u,v,w$ induce 
Fig.~\ref{Forbidden1}(d).
Hence $uw$ is a non-symmetric arc.
We claim that either $u$ or $w$ is a type 1 vertex and we prove this by contradiction.
Suppose that neither $u$ nor $w$ are type 1 vertices.
Assume first that $w$ is not a semi-strict di-simplicial vertex of $S(D)$, then $u$ must be a semi-strict di-simplicial vertex of $S(D)$.
Since $u$ is not a type 1 vertex, $u$ is type 2 and there exists a violating triple $(u_1,u,u_2)$ for $u$ such that exactly one of $u_1u$ or $uu_2$ is a symmetric arc. 
Without loss of generality, assume that $uu_2$ is a symmetric arc.
(The other case can be prove similarly by reversing the arcs of $D$.)
Then $u_1u$ is a non-symmetric arc, and, $u_1$ and $u_2$ are joining by a non-symmetric arc.
Since $w$ and $u_2$ are asynchronous neighbours of $u$, they are adjacent.
If $u_2w$ is a symmetric arc, then $vu_2$ is a symmetric arc as otherwise $u,v,w,u_2$ induce Fig.~\ref{Forbidden1}(b), a contradiction.
Hence $v$ and $u_1$ are adjacent because they are the asynchronous neighbours of $u_2$.
If $v$ and $u_1$ are joining by a symmetric arc then $v$ is not a semi-strict di-simplicial vertex of $S(D)$;
if $v$ and $u_1$ are joining by a non-symmetric arc then $(u_1,v,u_2)$ is a violating triple of $v$ which against the assumption that $v$ is a type 1 vertex.
Therefore, the arc between $u_2$ and $w$ is a non-symmetric arc.
By the assumption, $w$ is not a semi-strict di-simplicial vertex of $S(D)$, then there exist $w_1 \neq u_2,w_2 \neq u_2$ such that both $ww_1$ and $ww_2$ are symmetric arcs but $w_1$ and $w_2$ are not joining by a symmetric arc.
We can see that $u$ is adjacent to both $w_1$ and $w_2$.
Clearly, the arcs between $u$ and $w_1$ and between $u$ and $w_2$ are symmetric arcs, as otherwise $u,w,u_2,w_1$ or $u,w,u_2,w_2$ induce Fig.~\ref{Forbidden1}(a).
Since $w_1$ and $w_2$ are not adjacent by symmetric arc but both $uw_1$ and $uw_2$ are symmetric arcs, $u$ is not a semi-strict di-simplicial vertex of $S(D)$, which is a contradiction.
Therefore, if $w$ is not a semi-strict di-simplicial vertex, then $u$ is a type 1 vertex.
Due to the symmetric, if $u$ is not a semi-strict di-simplicial vertex, then $w$ is a type 1 vertex.
Thus, assume that both $u$ and $w$ are semi-strict di-simplicial vertices but none of them is type 1.
Then there exist violating triples $(u_1,u,u_2)$ for $u$ and $(w_1,w,w_2)$ for $w$.
Suppose that both $u_1u$ and $w_1w$ are non-symmetric arcs, then $uu_2$ are $ww_2$ are symmetric arcs.
(The other cases have the same proof idea due to the symmetric).
Again, we can see that $u_2$ is adjacent to both $v$ and $w$, and $w_2$ is adjacent to both $u$ and $v$.
By the previous discussion, $u_2$ and $w$ are joining by a non-symmetric arc.
Similarly, $w_2$ and $u$ are also joining by a non-symmetric arc.
Then the subdigraph induced by $u,w,u_2,w_2$ is Fig.~\ref{Forbidden1}(a), a contradiction.
Therefore, if $v$ is a type 1 vertex, then either $u$ or $w$ is also a type 1 vertex.

Now, suppose that $w$ is not a type 1 vertex.
Then $u$ must be a type 1 vertex, and there exists $u_1$ such that $(u_1,u,x_1)$ is a violating triple of $u$ and at lease one of the $u_1$ or $x_1$ is a type 1 vertex.
Since $v$ and $u_1$ are asynchronous neighbours of $u$, they are adjacent.
If they are joining by a symmetric arc, then $v$ is a type 2 vertex which contradicts with our assumption.
If $vu_1$ is a non-symmetric arc, then $u,v,u_1$ induce Fig.~\ref{Forbidden1}(d).
Therefore, $u_1v$ is a non-symmetric arc and so $u_1$ and $w$ are adjacent.
If $u_1w$ is a symmetric arc, then $u,v,w,u_1$ induce Fig.~\ref{Forbidden1}(c), a contradiction.
If $wu_1$ is a non-symmetric arc, then $v,w,u_1$ is Fig.~\ref{Forbidden1}(d), which is also a contradiction. 
Hence, $u_1w$ is a non-symmetric arc, and so $(u_1,v,w)$ is a violating triple for $v$. 
Since $w$ is not a type 1 vertex, $u_1$ is a type 1 vertex.
Thus, there exists a violating triple $(u_2,u_1,x_2)$ for $u_1$ where $u_2u_1$ and $u_1x_2$ are both non-symmetric arcs.
Similarly as above, we can see that both $u_2v$ and $u_2w$ are non-symmetric arcs and so $(u_2,v,w)$ is a violating triple for $v$.
Continuing this way in a finite number of steps, we end up with a directed cycle $u_1,u_2, \dots, u_k$, along with vertices $x_1, x_2, \dots, x_k$, such that for each $i = 1, 2, \dots, k$,
\begin{enumerate}
	\item $u_i$ is a semi-strict  di-simplicial vertex of $S(D)$ of type 1,
	\item $(u_i,v,w)$ is a violating triple, and
	\item $(u_{i+1},u_i,x_{i+1})$ is a violating triple where $u_{i+1}u_i$ is a non-symmetric arc (subscripts are modulo $k$).
\end{enumerate}
By Lemma \ref{noOrientedCycle} no such directed cycle consisting of non-symmetric arcs exists.
Hence, $w$ is a semi-strict di-simplicial vertex of $S(D)$ of type 1.
Let $(y_1,w,w_1)$ be a violating triple of $w$ where both $y_1w$ and $ww_1$ are non-symmetric arcs. 
Clearly, both $uw_1$ and $vw_1$ are non-symmetric arcs and so $(u,v,w_1)$ is a violating triple.
If $w_1$ is not a type 1 vertex, then we are back to the previous case. 
Hence, $w_1$ is a semi-strict di-simplicial vertex of $S(D)$ of type 1.
And so there exists a violating triple $(y_2,w_1,w_2)$ such that $(u,v,w_2)$ is a violating triple for $v$ and $w_2$ is a type 1 vertex.
Continuing this way, there exists a directed cycle $w_1,w_2, \dots, w_r$ along with vertices $y_1, y_2, \dots, y_r$, such that for each $i = 1, 2, \dots, r$,
\begin{enumerate}
	\item $w_i$ is a semi-strict di-simplicial vertex of $S(D)$ of type 1.
	\item $(u,v,w_i)$ is a violating triple, and 
	\item $(y_{i+1},w_i,w_{i+1})$ is a violating triple where $w_iw_{i+1}$ is a non-symmetric arc (subscripts are modulo $r$).
\end{enumerate}
Again, by Lemma \ref{noOrientedCycle}, no such directed cycle consisting of non-symmetric arcs exists. 
Therefore, any semi-strict di-simplicial vertex $v$ of $S(D)$ is a type 2 vertex.

Let $v$ be a semi-strict di-simplicial vertex of $S(D)$. 
Since $v$ is type 2, there is a violating triple $(u,v,w)$ such that exactly one of $uv$ or $vw$ is a non-symmetric arc. 
If $uv$ is a non-symmetric arc, then by Lemma \ref{canonical_WQTD} and the result that we just proved, $u$ is a semi-strict di-simplicial vertex of $S(D)$ of type 2.
If $vw$ is a non-symmetric arc, then $w$ is a semi-strict di-simplicial vertex of $S(D)$ of type 2.
Therefore, for any type 2 vertex $v$, there is a type 2 vertex $z$ such that $z$ and $v$ are part of a violating triple for $v$ and the arc between $z$ and $v$ is non-symmetric. 
It follows that there exists a circuit $z_1, z_2, \dots, z_l$ along with vertices $z_1',z_2', \dots, z_l'$, such that for each $i =  1, 2, \dots, l$, 
\begin{enumerate}
	\item $z_i$ is a semi-strict di-simplicial vertex of $S(D)$ of type 2.
	\item either $(z_{i+1},z_i,z_i')$ or $(z_i',z_i,z_{i+1})$ is a violating triple.
	\item the arc between $z_i$ and $z_{i+1}$ is non-symmetric, and the arc between  $z_i$ and $z_i'$ is symmetric (subscripts are modulo $l$).
\end{enumerate}
We assume that the circuit is chosen to have the minimum length. 
By Lemma \ref{noOrientedCycle} this circuit is not a directed cycle.
Assume without loss of generality that $z_1z_2$ and $z_1z_l$ are non-symmetric arcs.
Hence $(z_1',z_1,z_2)$ and $(z_1,z_l,z_l')$ are violating triples and the arcs between $z_2, z_1'$ and between $z_1, z_l'$ are non-symmetric.
Since $z_l$ and $z_1'$ are asynchronous neighbours of $z_1$, they are adjacent.
If they are adjacent by a non-symmetric arc, then $z_1,z_l',z_1',z_l$  induce 
Fig.~\ref{Forbidden1}(a), which is a contradiction.
Hence, $z_1'z_l$ is a symmetric arc, then $z_2$ and $z_l$ become asynchronous neighbours of $z_1'$ and so they are adjacent.
Also, $z_2z_l$ is not a symmetric arc as otherwise $z_1,z_2,z_l,z_1'$ induce 
Fig.~\ref{Forbidden1}(a).
If $z_2z_l$ is a non-symmetric arc, then $(z_2,z_l,z_1')$ is a violating triple and we end up with a smaller cycle, which contradicts with the assumption.
If $z_lz_2$ is a non-symmetric arc, then $(z_1',z_l,z_2)$ is a violating triple which is again a smaller cycle, a contradiction.
Therefore, no such circuit exists. 
$D$ has a semi-strict di-simplicial vertex. \qed
\end{proof}

\begin{corollary}
Let $D$ be an extended semicomplete digraph or a quasi-transitive digraph.
Then $D$ is semi-strict chordal if and only if $S(D)$ is semi-strict chordal and $D$ does not contain any digraph in Fig.~\ref{Forbidden1} as an induced subdigraph.
	\qed 
\end{corollary}

\section{Locally semicomplete semi-strict chordal digraphs}
\label{four}

In this section we characterize locally semicomplete semi-strict chordal digraphs
by forbidden subdigraphs.

\begin{lemma}\label{LSCDnecessity}
If $D$ is a semi-strict chordal digraph, then $D$ does not contain 
any digraph in Fig.~\ref{lollipop} as an induced subdigraph.
\end{lemma}

\begin{center} 
	\begin{figure}[htb]
		\center
		\begin{tikzpicture}[>=latex]
		\begin{pgfonlayer}{nodelayer}
		\node [style=vertex] (u1) at (0,0) {};
		\node [style=textbox] at (0, 0.5) {$u_1$};
		\node [style=vertex] (v1) at (-1.5, 1) {};
		\node [style=textbox] at (-2, 1) {$v_1$};
		\node [style=vertex] (v2) at (-1.5, -1) {};
		\node [style=textbox] at (-2, -1) {$v_2$};
		\node [style=vertex] (uk) at (3,0) {};
		\node [style=textbox] at (3, 0.5) {$u_k$};
		\node [style=vertex] (w1) at (4.5,1) {};
		\node [style=textbox] at (5, 1) {$w_1$};
		\node [style=vertex] (w2) at (4.5,-1) {};
		\node [style=textbox] at (5, -1) {$w_2$};
		\node [style=vertex] (u2) at (1,0) {};
		\node [style=textbox] at (1, 0.5) {$u_2$};
		\node [style=vertex] (u3) at (2,0) {};
		\node [style=textbox] at (2, 0.5) {$u_3$};
		\draw[style=arc] (v1) to (u1);
		\draw[style=arc] (v2) to (u1);
		\draw[style=arc] (v1) to (v2);
		\draw[style=arc] (v2) to (v1);
		\draw[style=arc] (uk) to (w1);
		\draw[style=arc] (uk) to (w2);
		\draw[style=arc] (w1) to (w2);
		\draw[style=arc] (w2) to (w1);
		\draw[style=arc] (u1) to (u2);
		\draw[style=arc] (u2) to (u3);
		\draw[style=dashed, ->] (u3) to (uk);
		
		\end{pgfonlayer}
		\end{tikzpicture}
		\caption{\label{lollipop} Locally semicomplete digraphs which are not semi-strict chordal}
	\end{figure}
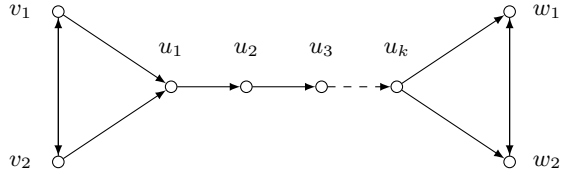
\end{center}
\begin{proof}
Suppose that a digraph in Fig.~\ref{lollipop} is an induced subdigraph in $D$.
Then the subdigraph of $D$ induced by the vertices of Fig.~\ref{lollipop} has no 
semi-strict di-simplicial vertex and hence is not a semi-strict chordal digraph, 
and so $D$ is not a semi-strict chordal digraph. Therefore, a semi-strict chordal 
digraph $D$ does not contain a digraph in Fig.~\ref{lollipop} as an induced 
subdigraph.
\qed
\end{proof}

\begin{lemma} \label{type1_LSCD}
Let $D$ be a locally semicomplete digraph such that $S(D)$ is semi-strict chordal 
and $D$ does not contain an directed cycle consisting of non-symmetric arcs, or 
any digraph in Fig.~\ref{Forbidden1} or Fig.~\ref{lollipop} as an induced 
subdigraph.
Suppose that $v$ is a semi-strict di-simplicial vertex of $S(D)$ and $(u,v,w)$ is 
a violating triple. Then the following statements hold:
\begin{enumerate}
\item If $v$ is of type 1, then both $u$ and $w$ are semi-strict di-simplicial 
vertices of $S(D)$ of type 1,
\item If $v$ is of type 2, then $u$ is a semi-strict di-simplicial vertex of $S(D)$
when $uv$ is a non-symmetric arc and $vw$ is a symmetric arc, and $w$ is a 
semi-strict di-simplicial vertex of $S(D)$ when $vw$ is a non-symmetric arc and 
$uv$ is a symmetric arc.
	\end{enumerate} 
\end{lemma}
\begin{proof}
We first prove statement~1: 
Since $v$ is of type 1, then both $uv$ and $vw$ are non-symmetric arcs.
If both $u$ and $w$ are semi-strict di-simplicial vertex of $S(D)$ of type 1, then 
we are done. So we may assume that at least one of $u, w$ is not of type 1.
Assume that $w$ is not a type 1 vertex. (A {\tiny }similar proof applies when $u$ is 
not a type 1 vertex.) Then $w$ is not a semi-strict di-simplicial vertex of S(D) or
is a type 2 vertex. In either case, there exist vertices $w_1, w_2$ such that $w_1w$ is 
a symmetric arc, $w_1$ and $w_2$ are joining by a non-symmetric arc. Clearly, $w_1 \neq v$. Since $v$ is an in-neighbour but not an 
out-neighbour of $w$ and $w_1$ is an out-neighbour of $w$, so $w_1$ and $v$ are 
adjacent. We consider two cases depending whether or not $u$ is a type 1 vertex. 

First suppose that $u$ is not a type 1 vertex. Then there exist vertices $u_1, u_2$
such that $uu_1$ is a symmetric arc, $u_1$ and $u_2$ are joining by a non-symmetric arc, $u$ and $u_2$ are adjacent, $u_1 \neq v$  and they are adjacent.
		We now suppose that $u$ and $w$ are not adjacent.
		Then $u_1$ is not adjacent to $w$, and $u$ is not adjacent to $w_1$, and so $u_1 \neq w_1$ and they are not adjacent.
		Moreover, both $u_1v$ and $vw_1$ are non-symmetric arcs.
		Then $u,v,w,u_1,w_1$ induce Fig.~\ref{lollipop}, a contradiction.
		Hence $uw$ is a non-symmetric arc.
		If $vu_1$ is a symmetric arc, then $v \neq u_2$ because $u_1$ and $u_2$ are joining by a non-symmetric arc.
		In addition, $v$ is joining with $u_2$.
		If $v$ and $u_2$ are joining by a non-symmetric arc, then $v$ is a type 2 vertex, which is a contradiction.
		If $vu_2$ is a symmetric arc, then $v$ is not a semi-strict di-simplicial vertex of $S(D)$, again a contradiction.
		Therefore, $v$ and $u_1$ are joining by a non-symmetric arc, and so $u_1w$ is a non-symmetric arc as otherwise $u,v,w,u_1$ induce Fig.~\ref{Forbidden1}(b) or (c).
		Similarly, we can see that $v$ and $w_1$ is joining by a non-symmetric arc, and so $uw_1$ is a non-symmetric arc.
		But then $u_1$ is adjacent to $w_1$ and $u,w,u_1,w_1$ induce Fig.~\ref{Forbidden1}(a), a contradiction.
		Therefore, for a type 1 vertex $v$ and any violating triple $(u,v,w)$, at lease one of the $u$ or $w$ is also type 1.

Suppose now that $u$ is a type 1 vertex,
		Then there exists $x_1$ such that $x_1u$ is a non-symmetric arc.
		We claim that $x_1$ is of type 1.
		If $u$ and $w$ are not adjacent, then $w_1$ is not joining to $u$.
		Hence, $vw_1$ is a non-symmetric arc, and the arc between $x_1$ and $w$ can not be symmetric. 
		Since $v$ is type 1, $x_1$ and $v$ are not joining by a symmetric arc.
		If $vx_1$ is a non-symmetric, then $u,v,x_1$ induce Fig.~\ref{Forbidden1}(d), a contradiction.
		If $x_1v$ is a non-symmetric arc, then since $(x_1,v,w)$ is a violating triple and $w$ is not a type 1 vertex, we have just proved that $x_1$ must be a type 1 vertex.
		If $x_1$ is not joining to $v$, then $x_1$ is not joining to $w$ as otherwise $u,v,w,x_1$ induce a directed cycle consisting of non-symmetric arcs.
		And so $x_1$ is not joining to $w_1$ for the same reason.
		Suppose that $x_1$ is not a type 1 vertex.
		Then there is a vertex $x_1'$ such that $x_1x_1'$ is a symmetric arc.
		Moreover, $x_1'$ is not joining to any of $v,w,w_1$, and so $x_1'u$ is a non-symmetric arc.
		But then the subdigraph induced by $u,v,w,x_1,x_1',w_1$ is Fig.~\ref{lollipop}, a contradiction.
		Hence, if $u$ and $w$ are not adjacent then $x_1$ is a type 1 vertex.
		Suppose that $uw$ is a non-symmetric arc.
		Then $u$ and $w_1$ are adjacent since $ww_1$ is a symmetric arc, 
		and again $vx_1$ is neither a non-symmetric arc nor a symmetric arc.
		If $x_1v$ is a non-symmetric arc, then $x_1w$ cannot be a symmetric arc as otherwise $u,v,w,x_1$ induce Fig.~\ref{Forbidden1}(c).
		Thus, as we proved above, $x_1$ is a type 1 vertex.
		If $x_1$ is not joining to $v$, then $x_1$ is not joining to $w$, and so not joining to $w_1$.
		Therefore, $uw_1$ is a non-symmetric arc.
		Assume that $x_1$ is not type 1.
		Then there exists $x_1'$ such that $x_1x_1'$ is a symmetric arc and $x_1'$ is not joining with any of $v,w,w_1$.	
		Hence, $x_1'u$ is a non-symmetric arc.
		But then $u,w,x_1,x_1',w_1$ induce Fig.~\ref{lollipop}, a contradiction.
		Therefore, $x_1$ is a type 1 vertex.
		Similarly, we can show that there exists a vertex $x_2$ 
		which is a semi-strict di-simplicial vertex of $S(D)$ of type 1, and $x_2x_1$ is a non-symmetric arc.
		Continuing this way in a finite number of steps, we will end up with a directed cycle $x_1,x_2, \dots, x_k$ along with vertices $y_1,y_2, \dots, y_k$ such that for each $i = 1, 2, \dots, k$,
		\begin{enumerate}
			\item $x_i$ is a semi-strict di-simplicial vertex of $S(D)$ of type 1,
			\item $(x_{i+1},x_i,y_i)$ is a violating triple, and
			\item $x_{i+1}x_i$ is a non-symmetric arc (subscripts are modulo $k$).
		\end{enumerate}
		Assume that this directed cycle consisting of non-symmetric arcs has the minimum length.
		Then there is a symmetric arc between a pair of non-consecutive vertices of the cycle. 
		Without loss of generality assume that $x_1x_s$ is a symmetric arc of the shortest distance along the cycle.
		Therefore, $x_2$ is adjacent to $x_s$. 
		Since $x_2$ and $x_s$ are adjacent, we have $s=3$.
		Then $(x_2,x_1,x_3)$ is a violating triple for $x_1$, which contradicts with the assumption that $x_1$ is of type 1.
		Hence, no such cycle exists. 
		Therefore, for any type 1 vertex $v$ with violating triple $(u,v,w)$, both $u$ and $w$ are type 1 vertices.

For statement~2, assume that $vw$ is a symmetric arc and $uv$ is a non-symmetric arc (the case for $uv$ is a symmetric arc can be proved similarly), and we are going to show that $u$ is a semi-strict di-simplicial vertex of $S(D)$.
Consider $S(D) - (N^-[w] \cap N^+[w])$. Since $u$ and $w$ are not joining by a symmetric arc, $u$ is not a vertex in $N^-[w] \cap N^+[w]$ and hence is a vertex in $S(D) - (N^-[w] \cap N^+[w])$.
Since $S(D)$ is semi-strict chordal, $S(D) - (N^-[w] \cap N^+[w])$ is also semi-strict chordal and so each component of $S(D) - (N^-[w] \cap N^+[w])$ contains a semi-strict di-simplicial vertex.
Let $u'$ be a semi-strict di-simplicial vertex of the component of $S(D) - (N^-[w] \cap N^+[w])$ that contains $u$.
If $u=u'$, then we are done.
Otherwise, there exists a directed path $u=u_1,u_2, \dots, u_k=u'$ in $S(D) - (N^-[w] \cap N^+[w])$.
Note that $w$ and $u_2$ are joining by a non-symmetric arc, and so $u,v,w,u_2$ induce Fig.~\ref{Forbidden1}(a), a contradiction.
Therefore no such directed path exists, and so $u$ is a semi-strict di-simplicial vertex of $S(D)$.
\qed
\end{proof}	

\begin{theorem}
	A locally semicomplete digraph $D$ is semi-strict chordal if and only if $S(D)$ is semi-strict chordal and $D$ does not contain as an induced subdigraph a directed cycle consisting of non-symmetric arcs or a digraph in Fig.~\ref{Forbidden1} or 
Fig.~\ref{lollipop}.
\end{theorem}
\begin{proof}
The necessity follows from Lemma \ref{dicycle}, Lemma \ref{necessity} and Lemma \ref{LSCDnecessity}.
For the sufficiency assume that $S(D)$ is semi-strict chordal and $D$ contains neither a directed cycle consisting of non-symmetric arcs nor a digraph in 
Fig.~\ref{Forbidden1} or Fig.~\ref{lollipop} as an induced subdigraph.
To prove $D$ is semi-strict chordal it suffices to show that $D$ has a semi-strict di-simplicial vertex.
Since $S(D)$ is chordal, $S(D)$ has semi-strict di-simplicial vertices.
If any of the semi-strict di-simplicial vertices of $S(D)$ is a semi-strict di-simplicial vertex of $D$ then we are done.
Hence, we assume that none of the semi-strict di-simplicial vertices of $S(D)$ is a semi-strict di-simplicial vertex of $D$.

First suppose that $S(D)$ has semi-strict di-simplicial vertices of type 1 and let $v$ be such vertex.
Then there exists a violating triple $(u,v,w)$ for $v$, and we know form Lemma \ref{type1_LSCD} that both $u$ and $w$ are semi-strict di-simplicial vertex of $S(D)$ of type 1.
Hence there exists a circuit $u_1,u_2, \dots u_l$ along with vertices $w_1,w_2, \dots, w_l$, such that for each $i = 1, 2, \dots, l$,
\begin{enumerate}
	\item $u_i$ is a semi-strict di-simplicial vertex of $S(D)$ of type 1,
	\item $(u_{i+1},u_i,w_{i+1})$ is a violating triple, and 
	\item $u_{i+1}u_i$ is a  non-symmetric arc (subscripts are modulo $l$).
\end{enumerate}
As the first part of the proof of Lemma \ref{type1_LSCD}, we know that no such circuit exists.
Therefore, every semi-strict di-simplicial vertex of $S(D)$ is of type 2.

Let $v$ be a semi-strict di-simplicial vertex of $S(D)$.
Since $v$ is of type 2, there is a violating triple $(u,v,w)$ such that exactly one of $uv$ or $vw$ is a non-symmetric arc.
By Lemma \ref{type1_LSCD}, if $uv$ is a non-symmetric arc, then $u$ is a semi-strict di-simplicial vertex of $S(D)$,
if $vw$ is a non-symmetric arc, then $w$ is a semi-strict di-simplicial vertex of $S(D)$.
This implies that for each semi-strict di-simplicial vertex $v$ of $S(D)$ there is a semi-strict di-simplicial vertex $x$ of $S(D)$ such that $x$ is part of a violating triple for $v$ and the arc between $x$ and $v$ is non-symmetric.
It follows that there exists a circuit $x_1,x_2, \dots, x_r$, along with $y_1,y_2,\dots, y_r$, such that for each $i= 1,2, \dots, r$,
\begin{enumerate}
	\item $x_i$ is a semi-strict di-simplicial vertex of $S(D)$ of type 2,
	\item either $(x_{i+1},x_i,y_i)$ or $(y_i,x_i,x_{i+1})$ is a violating triple, and 
	\item the arc between $x_i$ and $x_{i+1}$ is non-symmetric, and the arcs between $x_i$ and $y_i$ are symmetric (subscripts are modulo $r$).
\end{enumerate}
We again assume that the circuit is chosen to have the minimum length.
Suppose that $r=2$.
Then $x_1,x_2,y_1,y_2$ induce Fig.~\ref{Forbidden1}(a), which is a contradiction.
Hence, $r \geq 3$.
Clearly, $x_1$ and $y_2$ are adjacent.
If the arc between $x_1$ and $y_2$ is a non-symmetric arc, then $x_1,x_2,y_1,y_2$ induce Fig.~\ref{Forbidden1}(a).
If $x_1y_2$ is a symmetric arc, then $x_1$ is adjacent to $x_3$.
Since $x_1$ is a semi-strict di-simplicial vertex of $S(D)$, the arc between $x_1$ and $x_3$ cannot be symmetric arc.
Hence, $x_1$ and $x_3$ are joining by a non-symmetric arc, and so either$(x_3,x_1,y_2)$ or $(y_2,x_1,x_3)$ is a violating triple for $x_1$, which exists a circuit of length shorter than $r$, a contradiction.
Therefore, $D$ has a semi-strict di-simplicial vertex.
This completes the proof. \qed
\end{proof}

\section{Concluding remarks and open problems}
\label{five}

In this paper we characterized semi-strict chordal digraphs in terms of knotting
graphs. We also gave forbidden subdigraph characterizations of the digraphs within 
the class of locally semicomplete digraphs and the class of weakly quasi-transitive
digraphs. The forbidden subdigraphs for weakly quasi-transitive semi-strict chordal
digraphs are similar to the forbidden subdigraphs for weakly quasi-transitive
chordal digraphs \cite{hy}, and the forbidden subdigraphs for locally semicomplete
semi-strict chordal digraphs are similar to those for strict chordal digraphs. 
The latter result suggests there is a better chance to find a forbidden subdigraph
characterization for general semi-strict chordal digraphs. We propose as open 
problems for finding such a characterization as well as one for general chordal
digraphs.


\begin{thebibliography}{100}

\bibitem{bang} J. Bang-Jensen, Locally semicomplete digraphs: A generalization of
tournaments, J. Graph Theory 14 (1990) 371 - 390.

\bibitem{BH} J. Bang-Jensen and J. Huang, Quasi-transitive digraphs, 
J. Graph Theory 20 (1995) 141-161.

\bibitem{gallai} T. Gallai, Transitiv orientierbare graphen,
Acta Mathematica Academiae Scientiarum Hungarica 18 (1967) 25 - 66.

\bibitem{golumbic} M.C. Golumbic, {\em Algorithmic Graph Theory and Perfect Graphs}, Academic Press, New York (1980).

\bibitem{2h} P. Hell and C. Hern\'andez-Cruz, Strict chordal and strict split 
digraphs, Discrete Applied Math. 216 (2017) 609 - 617.

\bibitem{hy} J. Huang and Y.Y. Ye, Chordality of locally semicomplete and weakly 
quasi-transitive digraphs,
Discrete Mathematics 344 (2021) 112362.

\bibitem{kleitman} D.J. Kleitman, A note on perfect elimination digraphs,
SIAM J. Comput. 3 (1974) 280 - 282.

\bibitem{mckee} T.A. McKee, Strict chordal digraphs viewed as graphs with 
distinguished edges,
Discrete Applied Math. 247 (2018) 122 - 126.

\bibitem{meister} D. Meister and J.A. Telle, Chordal digraphs,
Theoret. Comput. Sci. 463 (2012) 73 - 83.

\bibitem{HR} L. Haskins and D.J. Rose, Toward characterization of perfect 
elimination digraphs, SIAM J. Comput. 2 (1973) 217 - 224.

\bibitem{rs} N. Robertson and P.D. Seymour, Graph minors III: Planar tree-width, 
J. Combinatorial Theory B 36 (1984) 49–64

\bibitem{ye} Y.Y. Ye, On chordal digraphs and semi-strict chordal digraphs,
Master thesis, University of Victoria, 2019.

\end{thebibliography}
\end{document}